\documentclass[11pt]{amsart}
\usepackage{epsfig}

\usepackage{graphicx}
\usepackage{amscd}
\usepackage{amsmath}
\usepackage{amsfonts}
\usepackage{yfonts}
\usepackage{psfrag}
\usepackage{amssymb}
\usepackage{verbatim}
\usepackage{comment}

\newtheorem{theorem}{Theorem}[section]
\theoremstyle{plain}

\newtheorem{corollary}[theorem]{Corollary}

\newtheorem{example}[theorem]{Example}

\newtheorem{lemma}[theorem]{Lemma}

\newtheorem{prop}[theorem]{Proposition}

\def\Hk{{\mathcal H}}
\def\Xk{\Xi}

\def\wtil{\widetilde}

\def\Pref{{\rm Pref}}

\newcommand{\gam}{\gamma}
\newcommand{\om}{\omega}
\def\Pmu{\P_{\!\!\mu}}
\def\Om{\Omega}
\newcommand{\Sig}{\Sigma}

\newcommand{\Gam}{\Gamma}

\newcommand{\R}{{\mathbb R}}

\def\N{{\mathbb N}}

\newcommand{\Prob}{{\mathbb P}\,}
\def\P{\Prob}
\newcommand{\Nat}{{\mathbb N}}

\def\bp{{\bf p}}

\def\Sk{{\mathcal S}}

\def\bx{{\mathbf x}}
\def\be{\begin{equation}}
\def\ee{\end{equation}}

\def\Gk{{\mathcal G}}

\newcommand{\es}{\emptyset}
\newcommand{\udim}{\overline{\dim}}
\def\ldim{\underline{\dim}}
\def\dimloc{{\dim}_{\rm loc}}

\def\ov{\overline}

\begin{document}

\title[Multiplicative fractals]{Hausdorff dimension for fractals invariant under the multiplicative integers}

\author{Richard Kenyon}
\address{Richard Kenyon\\ Department of Mathematics\\ Brown University\\ Providence, RI 02912}
\author{Yuval Peres}
\address{Yuval Peres\\ Microsoft Research}
\author{Boris Solomyak }
\address{Boris Solomyak, Box 354350, Department of Mathematics,
University of Washington, Seattle WA 98195}
\email{solomyak@math.washington.edu}

\begin{abstract}
We consider subsets of the (symbolic) sequence space that are invariant under the action of the semigroup of
multiplicative integers.
A representative example is the collection of all 0-1 sequences $(x_k)$ such that $x_k x_{2k}=0$ for all $k$. 
We compute the Hausdorff and Minkowski dimensions of these sets and show that they are typically different. The proof
proceeds via a variational principle for multiplicative subshifts.
\end{abstract}

\date{\today}

\thanks{
The research of R. K. 
was supported in part by NSF. The research of B. S. was partially supported by the NSF grant DMS-0968879.
}

\maketitle

\thispagestyle{empty}

\section{Introduction}

Central objects in symbolic dynamics and the theory of fractals are shifts of finite type, and more generally,
closed subsets of the symbolic space $\Sig_m := \{0,\ldots,m-1\}^\N$ that are invariant under  the shift
$
\sigma(x_1,x_2,x_3,\ldots)= (x_2,x_3,\ldots).
$
(we refer to them as ``subshifts'' for short).
To a subset $\Omega$ of $\Sig_m$ we can associate a subset of
$[0,1]$ by considering the collection of all reals whose base $m$ digit sequences belong to $\Omega$. Subshifts then
correspond to closed subsets of $[0,1]$ invariant under the map $x\mapsto mx$ (mod 1). It is known \cite{Furst} that all such
sets have the Hausdorff dimension equal to the Minkowski (box-counting) dimension, which is equal to $(\log m)^{-1}$ times
the topological entropy of $\sigma$ on $\Omega$.

Note that shift-invariance implies invariance under the action of the {\bf semigroup of additive positive integers.}
In contrast, in this paper
we consider subsets of $\Sig_m$ and the corresponding fractals in $[0,1]$, which arise from the action of the
{\bf semigroup of multiplicative integers.} Namely, given a subset $\Omega\subset\Sigma_m$ and an integer $q\ge 2$, we let
\begin{equation} \label{def-Xom}
X_\Om =X_\Om^{(q)}:= \Bigl\{\om =
{(x_k)}_{k=1}^\infty \in \Sig_m:\ {(x_{iq^\ell})}_{\ell=0}^\infty \in \Om\ \ \mbox{for all}\ i,\ q\nmid i\Bigr\}
\end{equation}
and consider the corresponding subset of $[0,1]$:
\be \label{def-Xiom}
\Xi_\Om:= \Bigl\{x = \sum_{k=1}^\infty x_k m^{-k}:\ {(x_k)}_1^\infty \in X_\Om\Bigr\}.
\ee
If $\Om$ is shift-invariant, then $X_\Om$ is invariant under the action of multiplicative integers:
$$
(x_k)_1^\infty \in X_\Om \ \Rightarrow\ {(x_{rk})}_{k=1}^\infty\in X_\Om \ \ \mbox{for all}\ r\in \N.
$$
If $\Om$ is a shift of finite type, we refer to $X_\Om$ (and $\Xi_\Om$) as the ``multiplicative shift of finite type.''

Our interest in these sets was prompted by work of Ai-Hua Fan, Lingmin Liao
and Jihua Ma \cite{Fan} 
who computed the Minkowski dimension  of the ``multiplicative golden mean shift''
\be \label{gold}
\Xi_g:= \Bigl\{x = \sum_{k=1}^\infty x_k 2^{-k}:\ x_k \in \{0,1\},\ x_k x_{2k}=0\ \ \mbox{for all } k\Bigr\}
\ee
and raised the question of computing
its Hausdorff dimension. They showed that the Minkowski dimension is
\be \label{goldM}
\dim_M(\Xi_g) = \sum_{k=1}^\infty \frac{\log_2 F_{k+1}}{2^{k+1}}= 0.82429\ldots,
\ee
where $F_k$ is the $k$-th Fibonacci number: $F_1=1,\ F_2 = 2, F_{k+1} = F_{k-1}+F_k$.
As a special case of our results we obtain the Hausdorff dimension $\dim_H(\Xi_g)$.
\begin{prop} \label{prop-gold}
We have
\be \label{gold1}
\dim_H(\Xi_g) = -\log_2 p = 0.81137\ldots,\ \ \mbox{where}\ p^3=(1-p)^2,\ \ \ 0<p<1
\ee
%$$p = \frac{1}{3}\Bigl[1 - 5{\Bigl(\frac{2}{11+3\sqrt{69}}\Bigr)}^{1/3} + {\Bigl(\frac{11+3\sqrt{69}}{2}\Bigr)}^{1/3} \Bigr]\,.$$
%

Thus, $\dim_H(\Xi_g) < \dim_M(\Xi_g).$
\end{prop}
%The set $\Xk_g$ (or rather, the set of its binary digit sequences which we denote $X_g$), is an example of what we call the
%``multiplicative shift of finite type.'' It is not shift invariant, but invariant under the action of
%multiplicative integers:
%$$
%(x_j)_1^\infty \in X_g \ \Rightarrow\ {(x_{rj})}_{j=1}^\infty\in X_g \ \ \mbox{for all}\ r\in \N.
%$$
%
%The usual sets defined in terms of their allowed digits (e.g.\ the classical middle-thirds Cantor set) all have their
%Minkowski and Hausdorff dimensions equal to each other \cite{Furst}.

Proposition~\ref{prop-gold} will follow from a more general result, Theorem~\ref{th-main} below. For an exposition which focuses on the set $\Xi_g$ see
\cite{KPS}.

\medskip

In order to visualize the set $\Xi_g$ we show the set $\wtil{\Xi}_g$ in Figure 1, which is obtained from $\Xi_g$ by the
transformation
$$
\sum_{k=1}^\infty x_k 2^{-k} \mapsto \Bigl( \sum_{k=1}^\infty x_{2k-1} 2^{-k}, \sum_{k=1}^\infty x_{2k} 2^{-k} \Bigr).
$$
It is easy to see that this transformation doubles the Minkowski and Hausdorff dimensions.
\begin{figure}[ht]
\includegraphics[height=2.4in]{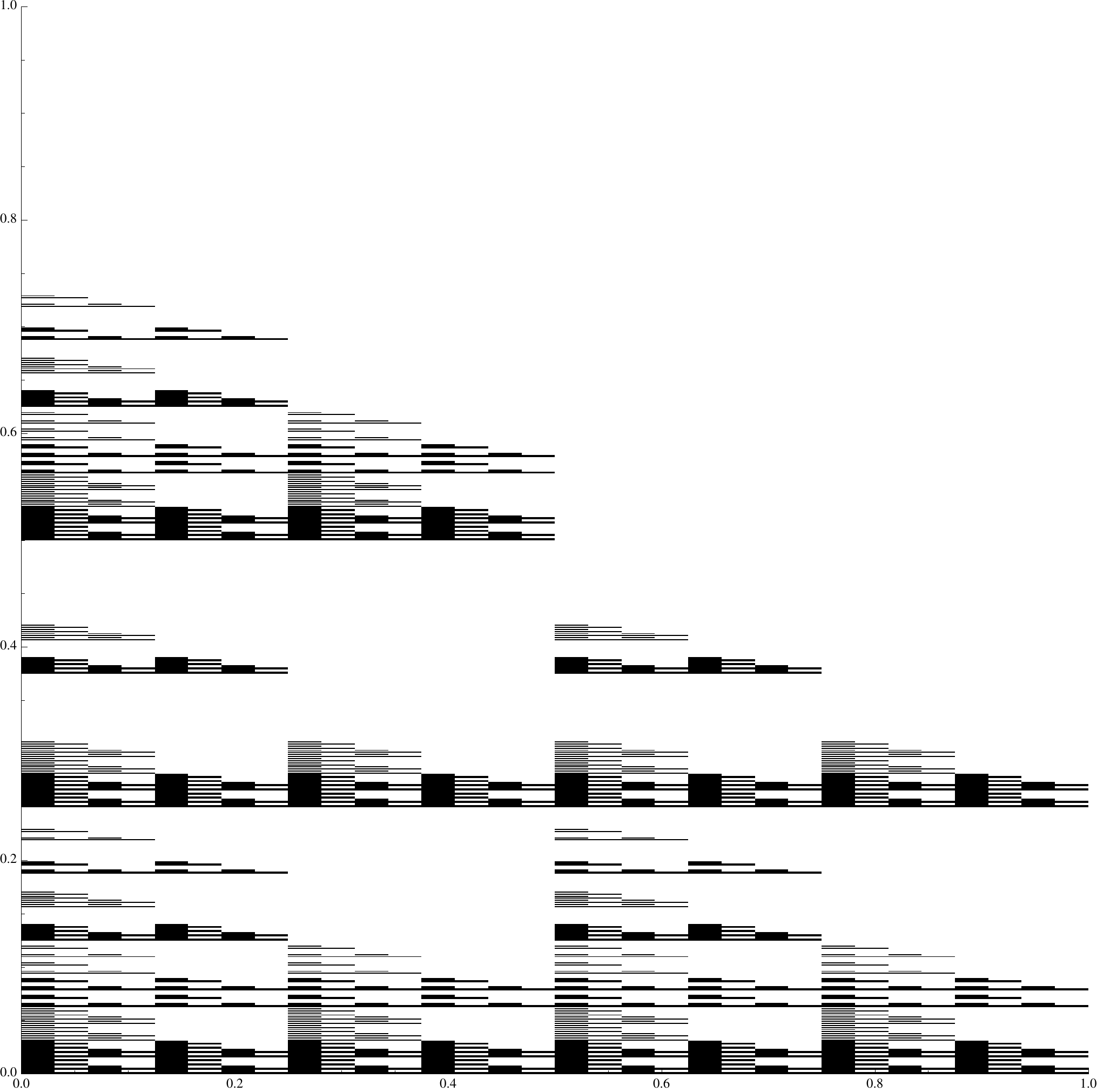}
\caption{Approximation of the set $\wtil{\Xk}_g$.}
\end{figure}

The figure resembles pictures of self-affine carpets, see \cite{Bedf,McM}, for which
the Hausdorff dimension is often less than the  Minkowski dimension. In fact, our proof bears some similarities with those of
\cite{Bedf,McM} as well.
An example of a self-affine set is shown in Figure 2.

\begin{figure}[ht]
\includegraphics[height=2.7in]{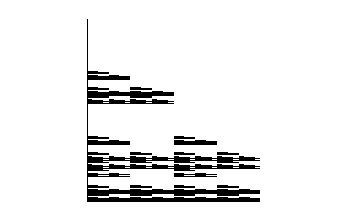}
\caption{Approximation of a self-affine set. Letting $M(x,y)=(\frac{x}2,\frac{y}4),$
 this set is $S=\left\{\sum_{i\ge0}M^i d_i :\, d_i
\in\{(0,0),(1,0),(0,1),(1,1),(0,3)\}\right\}.$}
\end{figure}

\medskip

The set $\Xi_g$ is a representative example of a large family of sets for which we compute the dimension.
Let $m\ge 2$ and let $A=(A(i,j))_{i,j=0}^{m-1}$ be a primitive (a non-negative matrix is primitive
if some power is strictly positive) matrix with 0-1 entries.
The usual (additive) shift of finite type determined by $A$ is defined as
$$
\Sig_A:= \{{(x_k)}_{k=1}^\infty:\ x_k\in \{0,\ldots,m-1\},\ A(x_k, x_{k+1}) =1,\ k\ge 1\}.
$$
Instead, we fix an integer $q\ge 2$ and consider the multiplicative shift of finite type 
\be \label{sft1}
X_A=X_A^{(q)}:= \{{(x_k)}_{k=1}^\infty:\ x_k\in \{0,\ldots,m-1\},\ A(x_k, x_{qk}) =1,\ k\ge 1\},
\ee
as well as the corresponding subset of the unit interval:
$$
\Xi_A:= \Bigl\{x = \sum_{k=1}^\infty x_k m^{-k},\ {(x_k)}_{k=1}^\infty\in X_A\}.
$$
%Observe that we get $\Xi_g$ by taking $q=2,\,m=2$, and $A = \left[\begin{array}{cc} 1 & 1 \\ 1 & 0 \end{array} \right]$.

As is well-known, the dimensions of $\Xi_A$ and $X_A$ coincide, if we use the standard metric on the
sequence space $\Sig_m$:
$$
\varrho\bigl((x_k), (y_k)\bigr) = m^{-\min\{n:\ x_n \ne y_n\}+1}
$$
on the
sequence space $\Sig_m$;
this is equivalent to restricting the covers of $\Xi_A$ to those by $m$-adic intervals. Thus, in the rest of
the paper we focus on the sets $X_A$.

In order to state our dimension result, we need the following elementary lemma.

\begin{lemma} \label{lem-elem}
Let $A=(A(i,j))_{i,j=0}^{m-1}$ be a primitive matrix, and $q>1$. Then there exists a unique vector $(t_i)_{i=0}^{m-1}$ satisfying
\be \label{eq-tvec}
t_i^q= \sum_{j=0}^{m-1} A(i,j) t_j,\ \ t_i>1,\ i=0,\ldots,m-1.
\ee
\end{lemma}

%Computing the Minkowski (box) dimension $\dim_M$ of $\Xi_A$ is straightforward.
 Below we use logarithms to base $m$, denoted $\log_m$, and write
$\ov{1}$ for the vector $(1,\ldots,1)^T \in \R^m$.

\begin{theorem} \label{th-main} {\rm (i)} Let $A$ be a primitive 0-1 matrix. Then the set $X_A$ given by (\ref{sft1}) satisfies
\be \label{eq-hausd1}
\dim_H(X_A) = \frac{q-1}{q} \log_m \sum_{i=0}^{m-1} t_i,
\ee
where $(t_i)_{i=0}^{m-1}$ is from Lemma~\ref{lem-elem}.

{\rm (ii)} The Minkowski dimension of $X_A$ exists and equals
\be \label{eq-Mink1}
\dim_M(X_A) = (q-1)^2 \sum_{k=1}^\infty \frac{\log_m(A^{k-1} \ov{1},\ov{1})}{q^{k+1}}\,.
\ee
We have $\dim_H(X_A) = \dim_M(X_A)$ if and only if $A$ has $\ov{1}$ as an eigenvector (i.e.\ row sums of $A$ are all equal).
\end{theorem}

The formula for the Minkowski dimension is not difficult to prove; it is included for comparison.

\subsection{Variational principle for multiplicative subshifts.}
We obtain Theorem~\ref{th-main} as a special case of a more general result. Let $\Om$ be an arbitrary
closed subset of $\Sig_m$ (it does not have to be shift-invariant), and define the sets $X_\Om$ and $\Xi_\Om$
by (\ref{def-Xom}) and (\ref{def-Xiom}). We refer to $X_\Om$ as a ``multiplicative subshift.''
Precise statements are given in the next section; here we just describe the results.

We can view our set $X_\Om$ as an infinite union of copies of $\Om$,
starting at all positive integers $i$ not divisible by $q$ (denoted $q\nmid i$) and ``sitting'' along
geometric progressions of ratio $q$. More precisely, denote $J_i = \{q^r i\}_{r=0}^\infty$ for $q\nmid i$ and let
$x|J_i = (x_{q^ri})_{r=0}^{\infty}$.
By definition (\ref{def-Xom}),
\be \label{cond2}
x\in X_\Om \ \Longleftrightarrow\ x|J_i \in \Om\ \ \mbox{for all}\ i,\ q\nmid i.
\ee
In order to compute (or estimate) the Hausdorff dimension of a set, one usually has to equip it with
a ``good'' measure and
calculate the appropriate ``H\"older exponent''. For subshifts, ``good'' measures are ergodic invariant measures.
For multiplicative subshifts,
their role is played by measures obtained in the following construction, essentially as an
infinite product of copies of a measure on $\Om$. Given a probability measure $\mu$ on $\Om$ we set
\be\label{eq-proba}
\Pmu[u] :=\prod_{i\le |u|,\, q\nmid i} \mu[u|J_i],
\ee
where $[u]$ denotes the cylinder set of all sequences starting with $u$ and 
$$
u|J_i = u_i u_{qi}\ldots u_{q^ri},\ \ q^ri\le |u| < q^{r+1}i.
$$
It is easy to verify that $\Pmu$ is a Borel probability measure supported on $X_\Om$ (see the next section 
for details).

For a probability measure $\P$, its Hausdorff dimension is defined by
$$
\dim_H(\P)=\inf\{\dim_H (F):\ F\ \mbox{Borel},\ \P(F)=1\},
$$
and the pointwise dimension at $x$ is given by
$$
\dimloc(\P,x) = \lim_{r\to 0} \frac{\log\P(B_r(x))}{\log r}\,,
$$
whenever the limit exists, where $B_r(x)$ denotes the open ball of radius $r$ centered at $x$.
We consider measures on the sequence space $\Sig_m$; then
\be \label{eq-ldim}
\dimloc(\P,x) = \lim_{n\to \infty} \frac{-\log\P[x_1^n]}{\log n},
\ee
where $x_1^n=x_1\ldots x_n$ denotes the initial segment (prefix) of the sequence $x$.
We prove that for any measure $\Pmu$ defined above, the pointwise dimension exists and
is constant $\Pmu$-a.e., which is then equal to $\dim_H(\Pmu)$ (see Proposition~\ref{prop-ldim}).
This can be viewed as a multiplicative analog of the Shannon-McMillan-Breiman Theorem and the 
entropy formula for the dimension of an ergodic shift-invariant measure $\nu$, namely,
$\dim_H (\nu)=h(\nu)/\log m$ (see \cite{Billing}).
Further, we obtain the ``Variational Principle for multiplicative subshifts,''
see Proposition~\ref{th-var}. 
%Then, finding the ``optimal measure'' yields the main result, that is, computing the Hausdorff dimension of
%$X_\Om$.
We can summarize this discussion with the following dictionary between the classical and multiplicative subshifts:

\bigskip

\begin{center} \begin{tabular}{|c|c|}
\hline {\bf classical} & {\bf multiplicative} \\ \hline \hline subshift $\Upsilon \subset \Sig_m$ & set $X_\Om$ \\ \hline
       invariant ergodic measure $\nu$ on $\Upsilon$ & measure $\Pmu$ \\ \hline
       Shannon-McMillan-Breiman Theorem & pointwise dimension of $\Pmu$  \\ \hline
       $\dim_H (\nu)=h(\nu)/\log m$ & dimension of $\Pmu$  \\ \hline
       Variational Principle: & $\dim_H(X_\Om) = \sup\{\dim_H(\Pmu):$ \\ 
       $\dim_H (\Upsilon) = \sup\{\dim_H(\nu):\ \nu$ is ergodic on $\Upsilon\}$ & \hfill $\mu$ is a probability on $\Om\}$\\ \hline
\end{tabular}
\end{center}

\bigskip

\section{General result. Variational problem.}

Let $\Om$ be an arbitrary
closed subset of $\Sig_m$,  and define the sets $X_\Om$ and $\Xi_\Om$
by (\ref{def-Xom}) and (\ref{def-Xiom}). Our general theorem computes the Hausdorff and Minkowski dimensions
of $X_\Om$ (as discussed earlier, the dimensions of $\Xi_\Om$ are the same as those of $X_\Om$).

Consider the {\bf tree of prefixes} of the set $\Om$. It is a directed graph $\Gam=\Gam(\Om)$ whose
set of vertices is
$$
V(\Gam) = \Pref(\Om) = \bigcup_{k=0}^\infty \Pref_k(\Om),
$$
where $\Pref_0(\Om)$ has only one element, the empty word $\varnothing$, and
$$\Pref_k(\Om) = \{u\in \{0,\ldots,m-1\}^k,\ \Om \cap [u]\ne \es\}.$$
There is a directed edge from a prefix $u$ to a prefix $v$ if $v=ui$ for some $i\in \{0,\ldots,m-1\}$.
In addition, there is an edge from $\varnothing$ to every $i\in \Pref_1(\Om)$.
Clearly, $\Gam(\Om)$ is a tree, and it has the outdegree bounded by $m$.
Note that if $\Om$ is shift-invariant, then the set $\Pref(\Om)$ coincides with the set of allowed (admissible)
words in $\Om$ (sometimes referred to as the {\bf language} of $\Om$).

The next lemma generalizes Lemma~\ref{lem-elem}.

\begin{lemma} \label{lem-elem2}
Let $\Gam = (V,E)$ be a directed graph (finite or infinite) with the outdegree bounded by $M<\infty$,
such that from each vertex there is at least one outgoing edge. Let $q>1$.
Then there exists a unique vector $\ov{t}\in [1,M^{\frac{1}{q-1}}]^V$ such that
\be \label{eq-tvec2}
t^q_v = \sum_{vw\in E} t_w,\ \ v\in V.
\ee
\end{lemma}

It is clear that Lemma~\ref{lem-elem} is a special case, with $\Gam$ being the directed graph with the
incidence matrix $A$.

Note that we only claim uniqueness of solutions in the given range.
In fact, uniqueness of positive solutions holds if we assume a priori bounds from zero and infinity; without
this assumption there may be infinitely many solutions on an infinite graph.

\begin{theorem} \label{th-main2} Let $\Om\subset \Sig$, and let $\ov{t}$ be the vector from Lemma~\ref{lem-elem2}
corresponding to the tree of prefixes $\Gam(\Om)$. Then

{\rm (i)}
\be\label{eq-hausd2}
\dim_H (X_\Om) = (q-1)\log_m t_\es;
\ee

{\rm (ii)}
\be \label{eq-mink2}
\dim_M (X_\Om) =
(q-1)^2 \sum_{k=1}^\infty \frac{\log_m|\Pref_k(\Om)|}{q^{k+1}}\,.
\ee
We have $\dim_H (X_\Om) = \dim_M (X_\Om)$ if and only if the
tree of prefixes is {\em spherically symmetric}, i.e.\ for every $k\in \Nat$, all prefixes of length $k$
have the same (equal) number of continuations in $\Pref_{k+1}(\Om)$.
\end{theorem}

Observe that Theorem~\ref{th-main} is a special case of  Theorem~\ref{th-main2}:
For part (i), we note that for a shift of finite type $\Sig_A$ the graph $\Gam(\Sig_A)$ has the
property that the tree of descendants of a prefix $u=u_1\ldots u_k$ depends only on the last symbol $u_k$.
Denote by $T_i$ this tree, which has $u_k=i$ as its root vertex, for $i=0,\ldots,m-1$,
and let $t_i$ be the solution of the system of equations (\ref{eq-tvec2}) evaluated at the root.
Here we use Lemma~\ref{lem-elem2}, with the uniqueness statement. Then we obtain from (\ref{eq-tvec2}) that
the vector ${(t_i)}_{i=0}^{m-1}$ satisfies (\ref{eq-tvec}). Finally, note that
$t^q_\varnothing = \sum_{i=0}^{m-1} t_i$ by (\ref{eq-tvec2}), hence
(\ref{eq-hausd2}) reduces to (\ref{eq-hausd1}).

For part (ii), we just note that $(A^{k-1} \ov{1},\ov{1})$ is the number of allowed words of length $k$
in the shift of finite type $\Sig_A$.

\subsection{Scheme of the proof. Statement of the Variational Principle.}
Recall (\ref{eq-proba}) that, given a probability measure $\mu$ on $\Om$ we define a measure on $X_\Om$ by
\be \label{eq-meas1}
\Pmu[u]:= \prod_{i\le n,\, q\nmid i} \mu[u|J_i], \ \ \mbox{where}\ |u|=n\ \mbox{and}\ J_i = \{q^r i\}_{r=1}^\infty.
\ee
This is a well-defined pre-measure on the semi-algebra of cylinder sets. Indeed, we have $\Pmu[i]=\mu[i]$ for
$i=0,\ldots,m-1$, and for $n+1=q^r i,\ q\nmid i$,
$$
\frac{\Pmu[u_1\ldots u_{n} u_{n+1}]}{\Pmu[u_1\ldots u_n]} = \frac{\mu[u_iu_{qi}\ldots u_{q^ri}]}{\mu[u_iu_{qi}\ldots u_{q^{r-1}i}]},
$$
whence
$$
\Pmu[u_1\ldots u_n] = \sum_{j=0}^{m-1} \Pmu[u_1\ldots u_{n}j].
$$
The extension of $\Pmu$ is a Borel measure supported on $X_\Om$,
since $\Om$ is a closed subset of $\Sig_m$ and
hence
$$
\Om = \bigcap_{k=1}^\infty \bigcup_{u \in \Pref_k(\Om)} [u].
$$

Observe that (\ref{eq-meas1}) is not the only way to put a measure on $X_\Om$: we could make the measure $\mu=\mu_i$
in (\ref{eq-meas1}) depend on $i$; however, this is not necessary for the purpose of computing the Hausdorff dimension.

We compute the Hausdorff dimension $\dim_H (\Pmu)$, which yields a lower bound on $\dim_H (X_\Om)$.
In order to state the result, we
need to introduce some notation.

For $k\ge 1$ let
$\alpha_k$ be the partition of $\Om$ into cylinders of length $k$:
$$
\alpha_k = \{\Om\cap [u]:\ u\in \Pref_k(\Om)\}=\{\Om \cap [u]:\ u\in \{0,\ldots,m-1\}^k,\ \Om \cap [u]\ne \es\}.
$$
For a measure $\mu$ on $\Sig_m$ and a finite partition $\alpha$, denote by $H_m^\mu(\alpha)$ the $\mu$-entropy of the
partition, with base $m$ logarithms:
$$
H_m^\mu(\alpha) = -\sum_{C\in \alpha} \mu(C)\log_m\mu(C).
$$
Now define
\be\label{def-smu}
s(\Om,\mu):= (q-1)^2 \sum_{k=1}^\infty \frac{H_m^\mu(\alpha_k)}{q^{k+1}}\,.
\ee

\begin{prop} \label{prop-ldim}
Let $\Om$ be a closed subset of $\Sig_m$ and $\mu$ a probability measure on $\Om$. Then
\be \label{eq-dimmeas}
\dimloc(\Pmu,x) = s(\Om,\mu)\ \ \mbox{for}\ \Pmu\mbox{-a.e.}\ x\in X_\Om.
\ee
Therefore,
$\dim_H(\Pmu) = s(\Om,\mu)$,  and $\dim_H(X_\Om) \ge s(\Om,\mu)$.
\end{prop}

We also have the Variational Principle:

\begin{prop}\label{th-var}
Let $\Om$ be a closed subset of $\Sig_m$. Then
\be \label{eq-var1}
\dim_H(X_\Om) = \sup_\mu \dim_H(\Pmu) = \sup_\mu s(\Om,\mu),
\ee
where the supremum is over Borel probability measures on $\Om$.
\end{prop}

It is clear from (\ref{def-smu}) that the function $\mu\mapsto s(\Om,\mu)$ is continuous on the compact 
space of probability measures with the
$w^*$-topology. Thus,
the supremum in (\ref{eq-var1}) is actually a maximum. Let
\be \label{def-som}
s(\Om) := \max\{s(\Om,\mu):\ \mu\ \mbox{is a probability on}\ \Om\}.
\ee
We call a measure $\mu$ for  which $s(\Om) = s(\Om,\mu)$
an {\bf optimal measure}. The next theorem characterizes such measures.

\begin{prop} \label{th-var2}
Let $\Om$ be a closed subset of $\Sig_m$ and let
$\ov{t}$ be the solution of the system of equations
(\ref{eq-tvec2}) for the tree of prefixes of $\Om$.
For any $k\ge 1$ and $u\in \Pref_k(\Om)$ let
\be \label{eq-meas12}
\mu[u] := \prod_{j=1}^k \frac{t_{u_1\ldots u_j}}{t^q_{u_1\ldots u_{j-1}}}\,.
\ee
This defines a probability measure $\mu$ on $\Om$. Moreover,

{\rm (i)} $\mu$ is the unique optimal measure;

{\rm (ii)} $s(\Om,\mu) = (q-1) \log_m t_\varnothing$.
\end{prop} 

Combining Propositions~\ref{th-var} and \ref{th-var2} yields part (i) of Theorem~\ref{th-main2}.

In the case when $\Om$ is a shift of finite type, the optimal measure turns out to be Markov.

\begin{corollary} \label{cor-var3}
Let $A$ be a primitive $m\times m$ 0-1 matrix and $\Sig_A$ the corresponding shift of finite type.
Let $\ov{t}={(t_i)}_{i=0}^{m-1}$ be the solution of the system of equations (\ref{eq-tvec}).
Then the unique optimal measure on
$\Sig_A$ is Markov, with the vector of initial probabilities $\bp = (\sum_{i=0}^{m-1} t_i)^{-1}\ov{t}$ and the matrix of transition
probabilities
$$(p_{ij})_{i,j=0}^{m-1}\ \mbox{where}\ p_{ij} = \frac{t_j}{t_i^q}\ \ \mbox{if}\ A(i,j)=1.$$
\end{corollary}

\section{Examples}

\begin{example}[golden mean] Let $q=2, m=2$, and $A = \left[\begin{array}{cc} 1 & 1 \\ 1 & 0 \end{array} \right]$.
Then $\Xk_A=\Xk_g$, the multiplicative golden mean shift from (\ref{gold}).
\end{example}

The system of equations (\ref{eq-tvec}) reduces to
$$
t_0^2=t_0+t_1,\ t_1^2=t_0,
$$
which immediately implies
$
t_1^3=t_1+1.
$
According to Corollary~\ref{cor-var3}, the optimal measure $\mu$ on $\Sig_A$ is Markov, with initial probability
of $0$ equal to $p = t_0/(t_0+t_1) = t_0^{-1}=t_1^{-2}$, and the initial probability of $1$ equal to
$1-p=t_1/(t_0+t_1) = t_1^{-3}$, whence $p^3=(1-p)^2$. The matrix of transition probabilities is
$\left[\begin{array}{cc} p & 1-p \\ 1 & 0 \end{array} \right]$.
Then, by
(\ref{eq-hausd1}),
$$
\dim_H(\Xk_g) = (1/2)\log_2(t_0+t_1)=-\log_2 p,
$$
which proves Proposition~\ref{prop-gold}. \qed

\begin{example}[Tribonacci] Let $q=2,\ m=3$, and $A = \left[\begin{array}{ccc} 1 & 1 & 1 \\ 1 & 0 & 0 \\ 0 & 1 & 0 \end{array} \right]$.
Then 
$$
X_A = \bigl\{{(x_k)}_1^\infty \in \{0,1,2\}^\Nat:\ x_k = 1\ \Rightarrow\ x_{2k}=0,\ x_k = 2 \ \Rightarrow\ x_{2k}= 1\bigr\}.
$$
We have
\be \label{eq-trib}
\dim_H(X_A) = 4\log_3 t \approx 0.726227,\ \ \mbox{where}\ t^4-t-1=0,
\ee
and
$$
\dim_M(X_A) = \sum_{k=1}^\infty \frac{\log_3 T_{k-1}}{2^{k+1}} \approx 0.75373,
$$
where $T_0 = 3,\ T_1 =5,\ T_2 =9,\ T_{k+2} = T_{k-1} + T_k + T_{k+1}$.
\end{example}

To verify (\ref{eq-trib}), we note that the equations (\ref{eq-tvec}) in this case are
$$
t_0^2 = t_0 + t_1 + t_2,\ t_1^2 = t_0,\ t_2^2 = t_1,
$$
whence $t_2^8 = t_2^4 + t_2^2 + t_2$. Thus, $t=t_2$ satisfies $t^7 = t^3+t+1$, and since $t^7-t^3-t-1 = (t^4-t-1)(t^3+1)$,
Theorem~\ref{th-main}(i) yields the formula for the Hausdorff dimension. The optimal measure is Markov, with the matrix
of transition probabilities equal to 
$\left[\begin{array}{ccc} t^{-4} & t^{-6} & t^{-7} \\ 1 & 0 & 0 \\ 0 & 1 & 0 \end{array} \right]$, and the vector of initial probabilities
$(t^{-4}, t^{-6}, t^{-7})$.

\begin{example}[2-step Markov] Let $q=2$, $m=2$, and 
$$
X:= \bigl\{{(x_k)}_1^\infty \in \{0,1\}^\N:\ x_k x_{2k} x_{4k} = 0,\ k\ge 1\bigr\}.
$$
Then $X = X_\Om$ where $\Om$ is the shift of finite type on the alphabet $\{0,1\}$ with the only forbidden 3-letter word $111$.
\end{example}

The graph  $\Gam(\Om)$ has the property that the tree of descendants of a prefix $u = u_1\ldots u_k$ of length $|u|\ge 2$
depends only on the last 
pair of symbols $u_{k-1}u_k=ij$. Denote by $T_{ij}$ this tree, for $i,j\in \{0,1\}$, and let $t_{ij}$ be the solution
of the system of equations (\ref{eq-tvec2}) evaluated at its root (we are using Lemma~\ref{lem-elem2} with uniqueness here). Then
(\ref{eq-tvec2}) on $\Gam(\Om)$ yields
\begin{eqnarray*}
t_{00}^2 & = & t_{00} + t_{01} \\ 
t_{01}^2 & = & t_{10} + t_{11} \\ 
t_{10}^2 & = & t_{00} + t_{01} \\
t_{11}^2 & = & t_{10},
\end{eqnarray*}
and also $t_\varnothing = t_{00}$.
Denoting $z=t_{11}$ we obtain, after a simple computation, that $(z^4-z^2)^2 = z^2 + z$ whence
$$
z^7 - 2z^5 + z^3 - z-1 = 0.
$$
Note also that $t_{00} = t_{11}^2 = z^2$.
Thus, by (\ref{eq-hausd2}),
$$
\dim_H(X) = 2\log_2 z \approx 0.956651.
$$

The Minkowski dimension in this example is 
$$\dim_M(X) = \sum_{k=1}^\infty \frac{\log_2R_{k-1}}{2^{k+1}}\approx 0.961789,$$
where $R_j$ is the number of allowed sequences of length $j$:
$R_1=2, R_2=4, R_3=7, R_{k+2}=R_{k-1}+R_k+R_{k+1}$.
\medskip

By the same method as in this example, one can
easily compute the Hausdorff dimension of $X_\Om$ where is an arbitrary (multi-step) shift of finite type.

\begin{sloppypar}
\begin{example}[Multiplicative $\beta$-shift] Let $\beta>1$ and $\Om = \Om_\beta$ be the $\beta$-shift (see \cite{Parry,DK} for the
definition and basic properties of $\beta$-shifts).
Let $q=2$. Then
\be \label{eq-beta}
\dim_H(X_{\Om_\beta}) = \log_2 t,\ \ \mbox{where}\ \ t = \sqrt{d_1t + \sqrt{d_2t + \sqrt{d_3t + \ldots}}}\,,
\ee
and $d_1 d_2 d_3\ldots$ is the infinite greedy expansion of $1$ in base $\beta$. Moreover, 
$\dim_M(X_{\Om_\beta})<\dim_H(X_{\Om_\beta})$
for all $\beta\not\in \N$.
\end{example}
\end{sloppypar}
Notice that $X_{\Om_\beta}$ is the 
multiplicative golden mean shift $X_g$ when $\beta = \frac{1+\sqrt{5}}{2}$, for which the infinite $\beta$-expansion of $1$ is
$101010\ldots$ 

The equation (\ref{eq-beta}) may be justified as follows.
Assume that $\beta\not\in \N$. % and let $m=\lfloor \beta \rfloor +1$. 
By \cite{Parry}, $x=
{(x_k)}_1^\infty\in \Om_\beta$
if and only if every shift of $x$, that is ${(x_k)}_\ell^\infty$ for $\ell\ge 2$, is less than or equal to ${(d_k)}_1^\infty$
in the lexicographic order. This implies that the tree of followers of the symbols $0,\ldots, d_1-1$ in $\Om_\beta$ is isomorphic to the
entire $\Pref(\Om_\beta)$, and we obtain the following equation at the root from (\ref{eq-tvec2}):
$$
t_\varnothing^2 = d_1 t_\varnothing + t_{d_1}.
$$
Similarly, we obtain 
$$
t_{d_1\ldots d_n}^2 = d_{n+1} t_\varnothing + t_{d_1\ldots d_{n+1}},\ \ n\ge 1,
$$
which easily reduces to (\ref{eq-beta}).

\section{Proof of Proposition~\ref{th-var2}}

Recall that for two partitions $\alpha$ and $\beta$, the conditional entropy is defined by
\be\label{conditionalentropydef}
H_m^\mu(\alpha|\beta) = \sum_{B\in \beta} \Bigl( - \sum_{A\in \alpha} \mu(A|B) \log_m \mu(A|B)\Bigr) \mu(B).
\ee

\noindent {\em Proof of Proposition~\ref{th-var2}(i).}
We have for $k\ge 2$,
\be \label{eq-condi1}
H_m^\mu(\alpha_k) = H_m^\mu(\alpha_k|\alpha_1) + H_m^\mu(\alpha_1),
\ee
by the properties of conditional entropy.
From (\ref{conditionalentropydef}),
$$
H_m^\mu(\alpha_k|\alpha_1) = \sum_{i=0}^{m-1} p_i H_m^{\mu_i}(\alpha_{k-1}(\Om_i)),
$$
where $p_i=\mu[i]$ and $H_m^{\mu_i}(\alpha_{k-1}(\Om_i))$ is the entropy of the partition of $\Om_i$,
the follower set of $i$ in $\Om$, into cylinders of length $k-1$, with respect to the measure $\mu_i$, which is the normalized
measure induced by $\mu$ on $\Om_i$.
Substituting this and (\ref{eq-condi1}) into (\ref{def-smu}) we obtain
\begin{eqnarray}
s(\Om,\mu) & = & \frac{q-1}{q} H_m^\mu(\alpha_1) + \frac{1}{q}\sum_{i=0}^{m-1} p_i s(\Om_i,\mu_i) \nonumber \\
           & = & \frac{q-1}{q} \Bigl[H_m^\mu(\alpha_1) + \frac{1}{q-1} \sum_{i=0}^{m-1} p_i s(\Om_i,\mu_i)\Bigr]. \label{eq-var3}
\end{eqnarray}
Now, the measure $\mu$ is completely determined by the probability vector $\bp=(p_i)_{i=0}^{m-1}$ and the conditional measures $\mu_i$. The
optimization problems on $\Om_i$ are independent, so if $\mu$ is optimal for $\Om$, then $\mu_i$ is optimal for $\Om_i$, for all
$i\le m$. Thus,
$$
s(\Om) = \max_\bp \frac{q-1}{q} \Bigl[H_m^\mu(\alpha_1) + \frac{1}{q-1} \sum_{i=0}^{m-1} p_i s(\Om_i)\Bigr].
$$
Observe that $H_m^\mu(\alpha_1) = -\sum_{i=0}^{m-1} p_i \log_m p_i$. It is well-known that
$$
\max_\bp \sum_{i=0}^{m-1} p_i(a_i-\log_m p_i) = \log_m\Bigl(\sum_{i=0}^{m-1} m^{a_i}\Bigr),
$$
which is achieved if and only if $p_i=m^{a_i}/\sum_{j=0}^{m-1} m^{a_j}$ for $i=0,\ldots,m-1$.
We have $a_i = s(\Om_i)/(q-1)$, which yields the optimal probability
vector
$$
\bp={(p_i})_{i=0}^{m-1},\ \ p_i = \frac{t_i}{t_\varnothing^q}\,,\ \ \mbox{where}
\ t_\varnothing:= m^{\frac{s(\Om)}{q-1}},\ \ t_i:=m^{\frac{s(\Om_i)}{q-1}},\ i\le m-1,
$$
and
$$
t_\varnothing^q = \sum_{i=0}^{m-1} t_i.
$$
This is the equation (\ref{eq-tvec2}) at the root of the graph $\Gam(\Om)$. However, the problem is analogous at each vertex, so
replacing the set $\Om$ with the set of followers of a prefix and repeating the argument, we obtain it for the entire graph.
We also get the formulas (\ref{eq-meas12}) for the optimal measure $\mu$ from the form of the optimal probability vector above.
Observe that the solution $\ov{t}$ of the system (\ref{eq-tvec2}) which we get this way is in the range $[1,m^{1/(q-1)}]$, where
we have uniqueness by Lemma~\ref{lem-elem2}. (Indeed, for any subtree $\Gam(\Om_u)$ of the tree $\Gam(\Om)$ we have the
outdegree bounded by $m$, and $s(\Om_u)\le 1$ by  (\ref{def-smu}) and (\ref{def-som}), in view of $H^\mu_m(\alpha_k)\le k$.)

This concludes the proof of Proposition~\ref{th-var2}(i), including the uniqueness statement.
\qed

\medskip

\noindent {\em Proof of Proposition~\ref{th-var2}(ii).}
In order to compute $s(\Om,\mu)$, it is useful to rewrite it in terms of conditional entropies. We have
$$
H_m^\mu(\alpha_{k+1}) = H_m^\mu(\alpha_k) + H_m^\mu(\alpha_{k+1}|\alpha_k).
$$
Applying this formula repeatedly, we obtain from (\ref{def-smu}):
\be \label{eq-dimsum}
s(\Om,\mu) = \sum_{k=1}^\infty \frac{(q-1)^2 H_m^\mu(\alpha_k)}{q^{k+1}} = \Bigl(\frac{q-1}{q}\Bigr)
\Bigl[H_m^\mu(\alpha_1) +\sum_{k=1}^\infty
\frac{H_m^\mu(\alpha_{k+1}|\alpha_k)}{q^k}\Bigr].
\ee
 Observe that
\begin{eqnarray*}
H_m^\mu(\alpha_1) & = & - \sum_{i=0}^{m-1} \frac{t_i}{t^q_\varnothing} \log_m\bigl( \frac{t_i}{t^q_\varnothing}\bigr) \\
                        & = & q \log_m t_\varnothing - \sum_{i=0}^{m-1} \frac{t_i}{t^q_\es} \log_m t_i
                          = q \log_m t_\varnothing - \sum_{i=0}^{m-1} \mu[i]\log_m t_i.
\end{eqnarray*}
Further,
\begin{eqnarray*}
H_m^\mu(\alpha_{k+1}|\alpha_k) & = & \sum_{[u]\in \alpha_k} \mu[u] \Bigl(- \sum_{j:\,[uj]\in \alpha_{k+1}}
\frac{t_{uj}}{t^q_u} \log_m\frac{t_{uj}}{t^q_u}\Bigr) \\
                                & = & \sum_{[u]\in \alpha_k} \mu[u] \Bigl(q\log_m t_u- \sum_{j:\,[uj]\in \alpha_{k+1}}
                                \frac{t_{uj}}{t^q_u} \log_m t_{uj}\Bigr)\\
                                & = & q \sum_{[u]\in \alpha_k} \mu[u] \log_m t_u - \sum_{[v]\in \alpha_{k+1}}
                                \mu[v]\log_m t_v,
\end{eqnarray*}
in view of $\mu[uj]= \mu[u]\frac{t_{uj}}{t^q_u}$. Now it is clear that the sum in (\ref{eq-dimsum}) telescopes, and
$s(\Om,\mu)=(q-1)\log_m t_\varnothing$, as desired. \qed

\medskip

We point out that Proposition~\ref{th-var2}(i) is not necessary for the proof of Theorem \ref{th-main2}, only 
Proposition~\ref{th-var2}(ii) is needed.

\section{Proof of the main theorem \ref{th-main2}}

\noindent {\em Proof of Proposition~\ref{prop-ldim}.}
Fix a probability measure $\mu$ on $\Om$. We are going to demonstrate that for every $\ell \in \N$, 
\be \label{eq-lb2}
\liminf_{n\to \infty} \frac{-\log_m \Pmu[x_1^n]}{n} \ge
(q-1)^2 \sum_{k=1}^\ell \frac{H_m^\mu(\alpha_k)}{q^{k+1}}\ \ \mbox{for $\Pmu$-a.e.}\ x,
\ee
and 
\be \label{eq-lb22}
\limsup_{n\to \infty} \frac{-\log_m \Pmu[x_1^n]}{n} \le
(q-1)^2 \sum_{k=1}^\ell \frac{H_m^\mu(\alpha_k)}{q^{k+1}}+ \frac{(\ell+1)\log_m(2m)}{q^\ell}\ \ \mbox{for $\Pmu$-a.e.}\ x.
\ee
Then, letting $\ell\to\infty$ will yield $\dimloc(\Pmu,x) = s(\Om,\mu)$ for $\Pmu$-a.e.\ $x$, as desired.

Fix $\ell\in \Nat$. To verify (\ref{eq-lb2}) and (\ref{eq-lb22}), we can restrict ourselves to $n=q^\ell r, \ r\in \N$.
(Indeed, if
$
q^\ell r \le n < q^\ell (r+1),
$
then
$$
\frac{-\log \Pmu[x_1^n]}{n} \ge
\frac{-\log \Pmu[x_1^{q^\ell r}]}{q^\ell(r+1)} \ge
\frac{r}{r+1} \cdot \frac{-\log \Pmu[x_1^{q^\ell r}]}{q^\ell r}\,,
$$
which implies that
$$
\liminf_{n\to \infty} \frac{-\log \Pmu[x_1^n]}{n} =
\liminf_{r\to \infty} \frac{-\log \Pmu[x_1^{q^\ell r}]}{q^\ell r}\,.
$$
The $\limsup$ is dealt with similarly.)

Let 
$$
\Gk_n = \Gk_{q^\ell r}:= \{j \le n:\ \exists\, i > n/q^\ell,\ q\nmid i,\ j \in J_i\}\ \ \mbox{and}\ \ \Hk_n:= \{j\le n:\ j\not\in \Gk_n\}.
$$
Then we have by the definition (\ref{eq-proba}) of the measure $\Pmu$:
\be \label{eq-lb33}
\Pmu[x_1^n] = \Pmu[x|\Gk_n]\cdot \Pmu[x|\Hk_n]
\ee
where $[x|\Gk_n]$ (resp.\ $[x|\Hk_n]$) denotes the cylinder set of $y\in X_\Om$ whose restriction to $\Gk_n$ (resp.\ $[x|\Hk_n]$)
coincides with that of $x$.

First we work with $\Pmu[x|\Gk_n]$. In view of (\ref{eq-proba}) we have
\be \label{eq-lb3}
\Pmu[x|\Gk_n] = \prod_{k=1}^\ell \prod_{\stackrel{\frac{n}{q^k} < i \le \frac{n}{q^{k-1}}}{q \nmid i}}
\mu[x_1^n|J_i].
\ee
Note that $x_1^n|J_i$ is a word of length $k$ for $i\in (n/q^k, n/q^{k-1}]$, $q\nmid i$, which is
a beginning of a sequence in $\Om$.
Thus, $[x_1^n|J_i]$ is an element of the partition $\alpha_k$. %Since $\Pmu$ is defined as a product measure,
The random variables $x\mapsto -\log_m\mu[x_1^n|J_i]$ are i.i.d\
for $i\in (n/q^k, n/q^{k-1}]$, $q\nmid i$,
and their expectation  equals $H_m^\mu(\alpha_k)$, by the definition of entropy.
Note that
\be \label{nondiv}
\#\bigl\{i\in (n/q^k, n/q^{k-1}]:\ q \nmid i\bigr\} =
\Bigl(\frac{q-1}{q}\Bigr)\Bigl(\frac{n}{q^k}-\frac{n}{q^{k-1}}\Bigr)= (q-1)^2\frac{n}{q^{k+1}}\,.
\ee
Fixing $k,\ell$ with $k\le \ell$ and taking $n=q^\ell r$,
$r\to \infty$, we get an infinite sequence of i.i.d.\ random variables. Therefore, by a version of the Law of Large
Numbers, we have
\be \label{eq-lb4}
\forall\ k\le \ell, \sum_{\stackrel{\frac{n}{q^k} < i \le \frac{n}{q^{k-1}}}{q \nmid i}}
\frac{-\log_m \mu[x_1^n|J_i]}{(q-1)^2(n/q^{k+1})} \
\to H_m^\mu(\alpha_k)\ \ \mbox{as}\ n = q^\ell r\to \infty,\ \ \mbox{for
$\Pmu$-a.e.\ $x$}.
\ee
By (\ref{eq-lb3}) and (\ref{eq-lb4}), for $\Pmu$-a.e.\ $x$,
\be \label{eq-lb44}
\frac{-\log_m\Pmu[x|\Gk_n]}{n} = \sum_{k=1}^\ell \frac{(q-1)^2}{q^{k+1}} \!\!\!\!\!
\sum_{\stackrel{\frac{n}{q^k} < i \le \frac{n}{q^{k-1}}}{q \nmid i}}
\frac{-\log_m\mu[x_1^n|J_i]}{(q-1)^2(n/q^{k+1})} \to \sum_{k=1}^\ell \frac{(q-1)^2H_m^\mu(\alpha_k)}{q^{k+1}}\,.
\ee
Since $\Pmu[x_1^n] \le \Pmu[x|\Gk_n]$, this proves (\ref{eq-lb2}).
Observe that (\ref{eq-lb2}) suffices for the lower bound $\dim_H(X_\Om) \ge \dim_H(\Pmu) \ge s(\Om,\mu)$, so the rest of the proof of this
proposition may be skipped if one is only interested in the computation of $\dim_H(X_\Om)$.

\medskip

Next we turn to (\ref{eq-lb22}), which requires working with $\Pmu[x|\Hk_n]$. In view of (\ref{nondiv}),
\begin{eqnarray}
|\Hk_n| = n - |\Gk_n|& = & n - \sum_{k=1}^\ell (q-1)^2 \frac{nk}{q^{k+1}}   \nonumber\\
                     & = & \frac{n}{q^\ell}\Bigl[(\ell+1) - \frac{\ell}{q}\Bigr] \label{exact}\\
                     &<& \frac{(\ell+1)n}{q^\ell}=(\ell+1)r.\label{lessthan} 
\end{eqnarray}
From (\ref{exact}),
\be \label{eq-BC}
\sum_{r=1}^\infty 2^{-|\Hk_{q^\ell r}|}< \infty.
\ee
Define
$$
\Sk(\Hk_n):= \bigl\{x\in X_\Om:\ \Pmu[x|\Hk_n] \le (2m)^{-|\Hk_n|} \bigr\}.
$$
Clearly, 
$$\Pmu(\Sk(\Hk_n)) \le 2^{-|\Hk_n|},
$$
since there are at most $m^{|\Hk_n|}$ cylinder sets $[x|\Hk_n]$. In view of (\ref{eq-BC}), 
$$\Pmu \bigl(\bigcap_{N\ge 1}\bigcup_{r=N}^\infty \Sk(\Hk_{q^\ell r})\bigr)=0,$$
hence
for $\Pmu$-a.e.\ $x\in X_\Om$ there exists $N(x)$ such that $x\not\in \Sk(\Hk_n)$ for all $n=q^\ell r\ge N(x)$.
For such $x$ and $n\ge N(x)$ we have (the last inequality from (\ref{lessthan}))
$$
\frac{-\log_m\Pmu[x|\Hk_n]}{n} < \frac{|\Hk_n|\log_m(2m)}{n} < \frac{(\ell+1)\log_m(2m)}{q^\ell}\,.
$$
Combining this with (\ref{eq-lb44}), which also holds $\Pmu$-a.e., and with (\ref{eq-lb33}), yields (\ref{eq-lb22}). \qed

\medskip
\noindent
{\em Proof of Proposition~\ref{th-var} and the upper bound in Theorem~\ref{th-main2}.} 
Often upper bounds for the Hausdorff dimension are obtained by explicit 
efficient coverings, which is easier than getting lower bounds.
This is not the case here, a feature shared with self-affine carpets from \cite{Bedf,McM}. In fact, we proceed similarly to
\cite{McM}, by exhibiting the ``optimal'' measure on the set $X_\Om$ to get an upper bound on the Hausdorff dimension.
We use the following well-known result; it essentially goes back to Billingsley \cite{Billing}.

\begin{prop}[see \cite{Falc}] \label{prop-mass}
Let $E$ be a Borel set in $\Sig_m$ and let $\nu$ be a finite Borel measure on $\Sig_m$.
If
$$\liminf_{n\to \infty} \frac{-\log_m \nu[x_1^n]}{n} \le s\ \ \mbox{for all}\ x\in E,$$
then $\dim_H(E) \le s$.
\end{prop}

It should be emphasized that the lower pointwise dimension of $\nu$ needs to be estimated from above {\bf for all} $x\in E$, unlike
in the proof of the lower bound, where the lower estimate for $\liminf$ is required only $\nu$-a.e.

\begin{lemma} \label{lem-average}
Let $\mu$ be the measure on $\Om$ defined by (\ref{eq-meas12}), and let $\Pmu$ be the
corresponding measure on $X_\Om$, defined by (\ref{eq-meas1}).
Then for any $x\in X_\Om$, denoting
$$
a_\ell(x):= \frac{-\log_m\Pmu[x_1^n]}{n}\ \ \mbox{for}\ n=q^\ell,
$$
we have
\be\label{eq-average}
\lim_{\ell\to \infty} \frac{a_1(x)+\cdots +a_\ell(x)}{\ell}=(q-1)\log_m t_\varnothing.
\ee
Thus, $\liminf_{\ell\to \infty} a_\ell(x) \le (q-1)\log_m t_\varnothing$ for all $x\in X_\Om$.
\end{lemma}

Once we prove the lemma, we are done with Theorem~\ref{th-main2}, since by Proposition~\ref{prop-mass} we will then get 
$\dim_H(X_\Om) \le (q-1)\log_mt_\varnothing$.  Proposition~\ref{th-var} then follows by Proposition~\ref{th-var2}(ii).

\medskip

\noindent{\em Proof of Lemma~\ref{lem-average}.}
Let $n=q^\ell$ and denote
$$
\bx_i^{(j)}:= x_i x_{qi} \cdots x_{q^ji}.
$$
We will also write $t(u)$ for $t_u$ in this proof, to make the formulas more readable.
Combining (\ref{eq-meas1}) with (\ref{eq-meas12}) yields
\begin{eqnarray}
-\log_m\Pmu[x_1^n] & = & - \sum_{k=1}^{\ell+1} \sum_{\stackrel{\frac{n}{q^k} < i \le \frac{n}{q^{k-1}}}{q \nmid i}}
\Bigl( \log_m \mu[x_i] + \sum_{j=1}^{k-1} \log_m
                  \frac{\mu[\bx_i^{(j)}]}{\mu[\bx_i^{(j-1)}]}\Bigr) \nonumber \\
                  & = & - \sum_{k=1}^{\ell+1} \sum_{\stackrel{\frac{n}{q^k} < i \le \frac{n}{q^{k-1}}}{q \nmid i}}
                  \Bigl( \log_m \frac{t(x_i)}{t^q(\varnothing)} + \sum_{j=1}^{k-1} \log_m
                  \frac{t(\bx_i^{(j)})}{t^q(\bx_i^{(j-1)})}\Bigr). \label{eq-ubd2}
\end{eqnarray}
For $\kappa \in \Nat$ and $x\in \Om$ denote
$$
\gam_x(\kappa) := \log_m t(\bx_i^{(j)}),\ \ \mbox{where}\ \kappa = q^ji,\ q\nmid i.
$$
Then, telescoping the sum $\sum_{j=1}^{k-1}$ in (\ref{eq-ubd2}) we obtain
$$
-\log_m\Pmu[x_1^n] = n(q-1)\log_m t(\varnothing) + (q-1)\sum_{\kappa=1}^{n/q} \gamma_x(\kappa) - \sum_{\kappa = n/q+1}^n \gam_x(\kappa).
$$
(Note that we pick up $q\log_m t(\varnothing)$
from each number in $[1,n]$ that is not divisible by $q$, for a total of $n(q-1)\log_m t(\varnothing)$.) Denote
$$
S_n:= \sum_{\kappa =1}^n \gam_x(\kappa);
$$
then
$$
-\log_m\Pmu[x_1^n] = n(q-1)\log_m t(\varnothing) + qS_{n/q}-S_n\ \ \mbox{for}\ n=q^\ell.
$$
We have for $n=q^\ell$, $\ell\ge 1$:
$$
a_\ell(x) = \frac{-\log_m\Pmu[x_1^n]}{n} = (q-1)\log_m t(\varnothing) + \frac{S_{n/q}}{n/q} - \frac{S_n}{n}\,.
$$
This implies
$$
\frac{a_1 + \cdots + a_\ell}{\ell} = (q-1)\log_m t(\varnothing) + \frac{S_1}{\ell} -  \frac{S_{q^\ell}}{\ell q^\ell}
\to (q-1)\log_m t(\varnothing),\ \ \mbox{as}\ \ell\to \infty,
$$
as desired.
\qed

\medskip

\noindent {\em Proof of Lemma~\ref{lem-elem2}.} We follow the scheme of the proof of \cite[Theorem 5.1]{Lyons}.

Let $V$ be the set of vertices of the graph and let $M$ be the maximal outdegree. Consider the space of functions
$Y:= [1,M^{1/(q-1)}]^V$ from $V$ to $[1,M^{1/(q-1)}]$,
which is compact in the topology of pointwise convergence, and the transformation $F:\,Y\to Y$, given by
$$
F(y_v) = \Bigl(\sum_{w:\ vw\in E} y_w\Bigr)^{1/q}.
$$
(It is easy to see that $F$ maps $Y$ into $Y$.)

Observe that $F$ is {\bf monotone} in the sense that
$$
\ov{y},\ov{z}\in Y,\ \ov{y}\le \ov{z}\ \Longrightarrow\ F(\ov{y}) \le F(\ov{z}),
$$
where ``$\le$'' is the pointwise partial order. Let $\ov{1}$ be the constant 1 function. Then
$\ov{1} \le F(\ov{1}) \le F^2(\ov{1})\le\ldots$ By compactness, there is a pointwise limit
$$
\ov{t} = \lim_{n\to \infty} F^n(\ov{1}),
$$
which is a fixed point of $F$, hence $\ov{t}$ satisfies the system of equations (\ref{eq-tvec2}).

It remains to verify uniqueness. Suppose $\ov{t}$ and $\ov{t'}$ are two distinct fixed points of $F$. Without
loss of generality, we can assume that $\ov{t}\not\le \ov{t'}$. Then let
$$
\alpha:= \inf\{\xi>1:\ \ov{t} \le \xi\ov{t'}\}.
$$
Clearly $\alpha\le M^{1/(q-1)}$. By continuity we have $\ov{t}\le \alpha\ov{t'}$, 
and so $1<\alpha$ by assumption. Now,
$$
\ov{t} = F(\ov{t}) \le F(\alpha\ov{t'}) = \alpha^{1/q} F(\ov{t'}) = \alpha^{1/q}\ov{t'},
$$
contradicting the definition of $\alpha$. The proof is complete. \qed

\medskip

\noindent {\em Proof of the statements on Minkowski dimension in Theorem~\ref{th-main2}.} It is well-known that one
can use covering by cylinder sets in the definition of lower Minkowski dimension, so we have for $X\subset \Sig_m$:
\be \label{lmink}
\ldim_M(X) = \liminf_{n\to \infty} \frac{\log_m|\Pref_n(X)|}{n}
\ee
where $|\Pref_n(X)|$ is the number of prefixes over all sequences in $X$; equivalently, the number of cylinder sets of length $n$ which intersect $X$.
We get the upper Minkowski dimension $\udim(X)$ by replacing $\liminf$ with $\limsup$ in (\ref{lmink}).

For dimension computations, we can restrict ourselves to $n$ from an arithmetic progression, so we can take
$n=q^\ell r$ for a fixed $\ell\in \Nat$.
Recall that $x\in X_\Om$ if and only if $x|J_i\in \Om$ for all $i$ such that $q \nmid i$. It follows that
$|\Pref_n(X_\Om)|$ is bounded below by the product of $|\Pref_k(\Om)|$ for each $i\in (n/q^k,n/q^{k-1}]$, with
$q\nmid i$, over $k=1,\ldots,\ell$. Thus, in view of (\ref{nondiv}), we have
$$
\log_m|\Pref_n(X_\Om)| \ge (q-1)^2 \sum_{k=1}^\ell \frac{n\log_m|\Pref_k(\Om)|}{q^{k+1}}\,.
$$
On the other hand,
$$
\log_m|\Pref_n(X_\Om)| \le (q-1)^2 \sum_{k=1}^\ell \frac{n\log_m|\Pref_k(\Om)|}{q^{k+1}} +
n - \sum_{k=1}^\ell k (q-1)^2 \frac{n}{q^{k+1}}
$$
by putting arbitrary digits in the remaining places. Dividing by $n$ and letting $n\to \infty$ we obtain
$$
\ldim_M(X_\Om) \ge (q-1)^2 \sum_{k=1}^\ell \frac{\log_m|\Pref_k(\Om)|}{q^{k+1}}
$$
and
$$
\udim_M(X_\Om) \le (q-1)^2 \sum_{k=1}^\ell \frac{\log_m|\Pref_k(\Om)|}{q^{k+1}} + (\ell+1)q^{-\ell}-
\ell q^{-\ell-1}.
$$
Since $\ell\in \Nat$ is arbitrary, this yields
(\ref{eq-mink2}). 

It remains to verify that $\dim_M(X_\Om) =  \dim_H(X_\Om)$ if and only if the tree of prefixes $\Gam(\Om)$ is
spherically symmetric. Compare the formula (\ref{eq-mink2}) with (\ref{def-smu}). Observe that 
$$
H_m^\mu(\alpha_k) \le \log_m|\Pref_k(\Om)|,
$$
with equality if and only every cylinder set $[u]$, for $u\in \Pref_k(\Om)$, has equal measure $\mu$.
To get $\dim_H(X_\Om)$, we have $\mu$ the optimal measure from (\ref{eq-meas12}). It is immediate from the equations
(\ref{eq-tvec2}) that the solution $t_u$ depends only on the length of the prefix $u$ if and only if $\Gam(\Om)$ is
spherically symmetric.
This implies the desired claim.
\qed

\section{Concluding remarks}

1. The motivation to consider the multiplicative golden mean shift $\Xi_g$ in \cite{Fan} came from the study of the dimension spectrum of certain multiple ergodic
averages. For $\theta\in [0,1]$ let 
$$
A_\theta = \Bigl\{{(x_k)}_1^\infty\in\Sig_2:\ \lim_{n\to\infty} \frac{1}{n} \sum_{k=1}^n x_k x_{2k} = \theta\Bigr\}.
$$
The authors of \cite{Fan} ask what is the Hausdorff dimension of $A_\theta$. It is easy to see that $\dim_H(A_0)=\dim_H(\Xi_g)$, and moreover, recently
 the methods developed
in the present paper have been adapted to compute the full dimension spectrum $\theta \mapsto \dim_H(A_\theta)$ \cite{PS}. Independently, the dimension of
$A_\theta$ and other sets of this type has been computed in \cite{FSW}.

\medskip

2. Not all subsets of $\Sig_m$ that are invariant under the action of multiplicative integers are of the form $X_\Om$ considered in this
paper. In fact, the sets of the form $X_\Om$ behave rather like full shifts,
because they are ``composed'' of independent copies of the set $\Om$, albeit in a ``staggered'' pattern.
On the other hand, let
$$
X:= \{x \in \Sig_2:\ x_k x_{2k} x_{3k} = 0\ \mbox{for all}\ k\}.
$$
Then clearly
$$
(x_k)_{k=1}^\infty \in X\ \Rightarrow\ {(x_{rk})}_{k=1}^\infty\in X \ \ \mbox{for all}\ r\in \N,
$$
but our methods are inadequate to compute the dimension of $X$.

\medskip

{\bf Acknowledgment.} We are grateful to J\"org Schmeling for passing to us the question about the Hausdorff dimension of the 
``multiplicative golden mean
shift'' $\Xi_g$.

\end{document}